\documentclass{article}

\newtheorem{theorem}{Theorem}[section]
\newtheorem{corollary}[theorem]{Corollary}

\newtheorem{lemma}[theorem]{Lemma}
\newtheorem{proposition}[theorem]{Proposition}
\newtheorem{definition}[theorem]{Definition}

\newtheorem{remark}{Remark}[section]

\evensidemargin\oddsidemargin
\begin{document}

\author{Vadim E. Levit and Eugen Mandrescu \\
%EndAName
Department of Computer Sciences\\
Holon Academic Institute of Technology\\
52 Golomb Str., P.O. Box 305\\
Holon 58102, ISRAEL\\
\{levitv, eugen\_m\}@barley.cteh.ac.il}
\title{On $\alpha $-Square-Stable Graphs}
\date{}
\maketitle

\begin{abstract}
The \textit{stability number }of a graph $G$, denoted by $\alpha (G)$, is
the cardinality of a maximum stable set, and $\mu (G)$ is the cardinality of
a maximum matching in $G$. If $\alpha (G)+\mu (G)$ equals its order, then $G$
is a K\"{o}nig-Egerv\'{a}ry graph. We call $G$ an $\alpha $-\textit{%
square-stable graph }if\textit{\ }$\alpha (G)=\alpha (G^{2})$, where $G^{2}$
denotes the second power of $G$. These graphs were first investigated by
Randerath and Wolkmann, \cite{ranvol}. In this paper we obtain several new
characterizations of $\alpha $-square-stable graphs\textit{. }We also show
that $G$ is an $\alpha $-square-stable K\"{o}nig-Egerv\'{a}ry graph if and
only if it has a perfect matching consisting of pendant edges. Moreover, we
find that well-covered trees are exactly $\alpha $-square-stable trees. To
verify this result we give a new proof of one Ravindra's theorem describing
well-covered trees, \cite{rav1}.
\end{abstract}

\section{Introduction}

All the graphs considered in this paper are simple, i.e., are finite,
undirected, loopless and without multiple edges. For such a graph $G=(V,E)$
we denote its vertex set by $V=V(G)$ and its edge set by $E=E(G).$ If $%
X\subset V$, then $G[X]$ is the subgraph of $G$ spanned by $X$. By $G-W$ we
mean the subgraph $G[V-W]$ , if $W\subset V(G)$. By $G-F$ we denote the
partial subgraph of $G$ obtained by deleting the edges of $F$, for $F\subset
E(G)$, and we use $G-e$, if $W$ $=\{e\}$. If $A,B$ $\subset V$ and $A\cap
B=\emptyset $, then $(A,B)$ stands for the set $\{e=ab:a\in A,b\in B,e\in
E\} $. A set of pairwise non-adjacent vertices is \textit{a stable set} of $%
G $. The \textit{stability number }of $G$, denoted by $\alpha (G)$, is the
cardinality of a maximum stable set in $G$. Let $\Omega (G)$ denotes $\{S:S$ 
\textit{is a maximum stable set of} $G\}$. $\theta (G)$ is the \textit{%
clique covering number} of $G$, i.e., the minimum number of cliques whose
union covers $V(G)$. Recall also that $i(G)=\min \{|S|:S\ $\textit{is a
maximal stable set in }$G\}$, and $\gamma (G)=\min \{|D|:D\ $\textit{is a
minimal domination set in }$G\}$. A \textit{matching} is a set of
non-incident edges of $G$. A matching of maximum cardinality $\mu (G)$ is a 
\textit{maximum matching}, and a \textit{perfect matching} is a matching
covering all the vertices of $G$. $M$ is an \textit{induced} \textit{matching%
}, \cite{cameron}, if no edge of $G$ connects two edges of $M$ (some recent
results on induced matchings can be found in \cite{gollask}, \cite{gollew}).
If $A,B$ are disjoint subsets of $V(G)$, we say that $A$ is \textit{uniquely
matched into} $B$ if there is a unique matching $M\subseteq (A,B)$ that
saturates all the vertices in $A$. $G$ is a \textit{K\"{o}nig-Egerv\'{a}ry
graph }provided $\alpha (G)+\mu (G)=\left| V(G)\right| $, \cite{dem}, \cite
{ster}. The neighborhood of a vertex $v\in V$ is the set $N(v)=\{w:w\in V$ \ 
\textit{and} $vw\in E\}$, and $N(A)=\cup \{N(v):v\in A\}$, for $A\subset V$.
If $G[N(v)]$ is a complete subgraph in $G$, then $v$ is a \textit{simplicial
vertex} of $G$. A maximal clique in $G$ is called a \textit{simplex} if it
contains at least one simplicial vertex of $G$, \cite{chhala}. $G$ is said
to be \textit{simplicial} if every vertex of $G$ is simplicial or is
adjacent to a simplicial vertex of $G$, \cite{chhala}. If $\left|
N(v)\right| =\left| \{w\}\right| =1$, then $v$ is a \textit{pendant vertex}
and $vw$ is a \textit{pendant edge} of $G$. By $C_{n}$, $K_{n}$, $P_{n}$ we
denote the chordless cycle on $n\geq $ $4$ vertices, the complete graph on $%
n\geq 1$ vertices, and respectively the chordless path on $n\geq 3$
vertices. A graph $G$ is $\alpha ^{-}$\textit{-stable} if $\alpha
(G-e)=\alpha (G)$, for any $e\in E(G)$, and $\alpha ^{+}$-\textit{stable} if 
$\alpha (G+e)=\alpha (G)$, for any edge $e\in E(\overline{G})$, where $%
\overline{G}$ is the complement of $G$, \cite{gun}. $G$ is \textit{%
well-covered} if it has no isolated vertices and if every maximal stable set
of $G$ is also a maximum stable set, i.e., it is in $\Omega (G)$, \cite{plum}%
. $G$ is called \textit{very well-covered}, \cite{fav1}, provided $G$ is
well-covered and $\left| V(G)\right| =2\alpha (G)$.

The distance between two vertices $v,w\in V(G)$ is denoted by $dist_{G}(v,w)$%
, or $dist(v,w)$ if no ambiguity. $G^{2}$ denotes the second power of graph $%
G$, i.e., the graph with the same vertex set $V$ and an edge is joining
distinct vertices $v,w\in V$ whenever $dist_{G}(v,w)\leq 2$. Clearly, any
stable set of $G^{2}$ is stable in $G$, as well, while the converse is not
generally true. Therefore, we may assert that $1\leq \alpha (G^{2})\leq
\alpha (G)$. Let notice that the both bounds are tight. For instance, if $G$
is not a complete graph and $dist(a,b)\leq 2$ holds for any $a,b\in V(G)$,
then $\alpha (G)\geq 2>1=\alpha (G^{2})$, e.g., for the $n$-star graph $%
G=K_{1,n}$, with $n\geq 2$, we have $\alpha (G)=n>$ $\alpha (G^{2})=1$. On
the other hand, if $G=P_{4}$, then $\alpha (G)=\alpha (G^{2})=2$.

In this paper we characterize the graphs $G$ for which the upper bound of
the above inequality is achieved, i.e., $\alpha (G)=\alpha (G^{2})$. These
graphs we call $\alpha $-\textit{square-stable,} or shortly square-stable.
We show that any square-stable graph is $\alpha ^{+}$-stable and that none
of them is $\alpha ^{-}$-stable. We give a complete description of
square-stable K\"{o}nig-Egerv\'{a}ry graphs extending the investigation of
well-covered trees, started in \cite{rav1}.

Randerath and Volkmann, \cite{ranvol}, prove that:

\begin{theorem}
\label{th3}\cite{ranvol} For a graph $G$ the following statements are
equivalent:

($\mathit{i}$) every vertex of $G$ belongs to exactly one simplex of $G;$

($\mathit{ii}$) $G$ satisfies $\alpha (G)=\alpha (G^{2})$;

($\mathit{iii}$) $G$ satisfies $\theta (G)=\theta (G^{2})$;

($\mathit{iv}$) $G$ satisfies $\alpha (G^{2})=\theta (G^{2})=\gamma
(G)=i(G)=\alpha (G)=\theta (G)$.
\end{theorem}

\begin{remark}
In general, it can be shown (e.g., see \cite{ranvol}) that the graph
invariants appearing in the above theorem are related by the following
inequalities: 
\[
\alpha (G^{2})\leq \theta (G^{2})\leq \gamma (G)\leq i(G)\leq \alpha (G)\leq
\theta (G).
\]

The graph $C_{12}$ indicates that no other non-trivial equality (except $%
\alpha (G)=\alpha (G^{2})$ and $\theta (G)=\theta (G^{2})$) of a pair of the
above invariants ensures that all of them are equal, namely, $\alpha
(C_{12}^{2})=i(C_{12})=4$, while $\alpha (C_{12})=6$.
\end{remark}

The following characterization of maximum stable sets in a graph, due to
Berge, we shall use in the sequel.

\begin{proposition}
\label{prop1}\cite{berge1} $S\in \Omega (G)$ if and only if every stable set 
$A$ of $G$, disjoint from $S$, can be matched into $S$.
\end{proposition}

Other useful results are:

\begin{proposition}
\label{prop8}\cite{levm} A graph $G$ is very well-covered if and only if it
is a well-covered K\"{o}nig-Egerv\'{a}ry graph.
\end{proposition}

\begin{proposition}
\label{prop12}\cite{levm} A K\"{o}nig-Egerv\'{a}ry graph is well-covered if
and only if it is very well-covered.
\end{proposition}

\begin{proposition}
\label{prop11}\cite{hayn} A graph $G$ is:

($\mathit{i}$) $\alpha ^{+}$-stable if and only if $\left| \cap \{S:S\in
\Omega (G)\}\right| \leq 1$;

($\mathit{ii}$) $\alpha ^{-}$-stable if and only if $\left| N(v)\cap
S\right| \geq 2$ is true for every $S\in \Omega (G)$ and any $v\in V(G)-S$.
\end{proposition}

By Proposition \ref{prop11}, an $\alpha ^{+}$-stable graph may have either $%
\left| \cap \{S:S\in \Omega (G)\}\right| =0$ or $\left| \cap \{S:S\in \Omega
(G)\}\right| =1$. This motivates the following definition.

\begin{definition}
\cite{levm} A graph $G$ is called:

($\mathit{i}$) $\alpha _{0}^{+}$-stable whenever $\left| \cap \{S:S\in
\Omega (G)\}\right| =0$;

($\mathit{ii}$) $\alpha _{1}^{+}$-stable provided $\left| \cap \{S:S\in
\Omega (G)\}\right| =1$.
\end{definition}

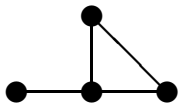
\begin{figure}[h]
\setlength{\unitlength}{1cm}%
\begin{picture}(5,1)\thicklines

  \multiput(6,0)(1,0){3}{\circle*{0.29}}
  \put(7,1){\circle*{0.29}}
  \put(6,0){\line(1,0){2}}
  \put(7,0){\line(0,1){1}} 
  \put(8,0){\line(-1,1){1}} 
  
 \end{picture}
\caption{Graph\emph{\ }$K_{3}+e.$}
\label{fig2}
\end{figure}

For instance, the graph in Figure \ref{fig2} is an $\alpha _{1}^{+}$-stable
graph.

In \cite{hayn} it was shown that an $\alpha ^{+}$-stable tree can be only $%
\alpha _{0}^{+}$-stable, and this is exactly the case of trees possessing a
perfect matching. This result was generalized to bipartite graphs in \cite
{levm1}. Nevertheless, there exist both $\alpha _{1}^{+}$-stable
K\"{o}nig-Egerv\'{a}ry graphs (e.g., the graph in Figure \ref{fig2}), and $%
\alpha _{0}^{+}$-stable K\"{o}nig-Egerv\'{a}ry graphs (e.g., all $\alpha
^{+} $-stable bipartite graphs).

\section{Square-stable graphs}

Clearly, any complete graph is square-stable. Moreover, since $%
K_{n}^{2}=K_{n}$, we get that 
\[
\Omega (K_{n})=\Omega (K_{n}^{2})=\{\{v\}:v\in V(K_{n})\}. 
\]

\begin{proposition}
\label{prop2}Graph $G$ is square-stable if and only if $\Omega
(G^{2})\subseteq \Omega (G)$.
\end{proposition}

\setlength {\parindent}{0.0cm}\textbf{Proof.} Clearly, any stable set $A$ of 
$G^{2}$ is stable in $G$, too. Consequently, if $G$ is square-stable, then
any maximum stable set of $G^{2}$ is a maximum stable set of $G$, as well,
i.e., $\Omega (G^{2})\subseteq \Omega (G)$.\setlength {\parindent}{3.45ex}

The converse is obvious. \rule{2mm}{2mm}\newline

It is quite evident that $G$ and $G^{2}$ are simultaneously connected or
disconnected. In addition, if $H_{i}$, $1\leq i\leq k$ are the connected
components of graph $G$, then $S\in \Omega (G)$ if and only if $S\cap
V(H_{i})\in \Omega (H_{i}),1\leq i\leq k$. Henceforth, using Proposition \ref
{prop2} we infer that:

\begin{proposition}
A disconnected graph is square-stable if and only if any of its connected
components is square-stable.
\end{proposition}

Therefore, in the rest of the paper all the graphs are connected, unless
otherwise stated.

\begin{lemma}
\label{lem4}Every $S\in \Omega (G^{2})$ has the property that 
\[
dist_{G}(a,b)\geq 3\ holds\ for\ any\ distinct\ a,b\in S.
\]
\end{lemma}

\setlength {\parindent}{0.0cm}\textbf{Proof.} If $S\in \Omega (G^{2})$ and $%
a,b\in S$, then $dist_{G}(a,b)\geq 3$, since otherwise $ab\in E(G^{2})$,
contradicting the stability of $S$ in $G^{2}$. \rule{2mm}{2mm}%
\setlength
{\parindent}{3.45ex}

\begin{proposition}
\label{prop3}A graph $G$ is square-stable if and only if there is some $S\in
\Omega (G)$ such that $dist_{G}(a,b)\geq 3$ holds for any distinct $a,b\in S$%
.
\end{proposition}

\setlength {\parindent}{0.0cm}\textbf{Proof.} If $G$ is square-stable, then
Proposition \ref{prop2} ensures that $\Omega (G^{2})\subseteq \Omega (G)$,
and by above Lemma \ref{lem4}, $dist(a,b)\geq 3$ is valid for every $S\in
\Omega (G^{2})$ and any $a,b\in S$.\setlength {\parindent}{3.45ex}

Conversely, if $S\in \Omega (G)$ and $dist(a,b)\geq 3$ holds for any $a,b\in
S$, then $S$ is stable in $G^{2}$, and therefore, $\left| S\right| \leq
\alpha (G^{2})\leq \alpha (G)=|S|$ implies that $\alpha (G^{2})=\alpha (G)$,
i.e., $G$ is square-stable. \rule{2mm}{2mm}

\begin{lemma}
\label{lem3}If $G\neq K_{\left| V\right| }$ is square-stable, then for every 
$S\in \Omega (G^{2})$ and any $a\in S$, there is $b\in S$ with $%
dist_{G}(a,b)=3$.
\end{lemma}

\setlength {\parindent}{0.0cm}\textbf{Proof.} Suppose, on the contrary, that
there are $S\in \Omega (G^{2})$ and some $a\in S$, such that $%
dist_{G}(a,b)\geq 4$ holds for any $b\in S$. Let $v\in V$ be with $%
dist_{G}(a,v)=2$; hence $dist_{G}(v,w)\geq 2$ is valid for any $w\in S$, and
consequently, $S\cup \{v\}$ is stable in $G$, a contradiction, because $S$
is a maximum stable set in $G$. \rule{2mm}{2mm}%
\setlength
{\parindent}{3.45ex}

\begin{lemma}
\label{lem5}If $G$ is square-stable, then $\Omega (G^{2})=\Omega (G)$ if and
only if $G$ is a complete graph.
\end{lemma}

\setlength {\parindent}{0.0cm}\textbf{Proof.} Suppose, on the contrary, that 
$\Omega (G^{2})=\Omega (G)$ holds for a non-complete square-stable graph $G$%
. Let $S\in \Omega (G^{2})$ and $a\in S$. According to Lemma \ref{lem3},
there is $b\in S$ with $dist_{G}(a,b)=3$. Now, if $c\in N(a)$ and $%
dist_{G}(c,b)=2$, Proposition \ref{prop3} implies that $S\cup \{c\}-\{a\}\in
\Omega (G)-\Omega (G^{2})$, contradicting the relation $\Omega
(G^{2})=\Omega (G)$.\setlength {\parindent}{3.45ex}

The converse is clear. \rule{2mm}{2mm}\newline

Combining Proposition \ref{prop2} and Lemma \ref{lem5} we obtain the
following assertion:

\begin{theorem}
\label{th4}$\Omega (G^{2})=\Omega (G)$ if and only if $G$ is a complete
graph.
\end{theorem}

Let $A\bigtriangleup B$ denotes the symmetric difference of the sets $A,B$,
i.e., the set 
\[
A\bigtriangleup B=(A-B)\cup (B-A). 
\]

\begin{theorem}
\label{th1}For a graph $G$ the following assertions are equivalent:

($\mathit{i}$) $G$ is square-stable;

($\mathit{ii}$) there exists $S_{0}\in \Omega (G)$ that satisfies the
property

$P1$: \textit{any stable set }$A$ \textit{of} $G$ \textit{disjoint from }$%
S_{0}$ \textit{can be uniquely matched into }$S_{0}$;

($\mathit{iii}$) any $S\in \Omega (G^{2})$ has property $P1$;

($\mathit{iv}$) for any $S_{1}\in \Omega (G)$ and $S_{2}\in \Omega (G^{2})$, 
$G[S_{1}\bigtriangleup S_{2}]$ has a unique perfect matching;

($\mathit{v}$) for any $S_{1}\in \Omega (G)$ and $S_{2}\in \Omega (G^{2})$, $%
G[S_{1}\bigtriangleup S_{2}]$ has a perfect matching;

($\mathit{vi}$) for any $S_{1}\in \Omega (G)$ and $S_{2}\in \Omega (G^{2})$, 
$G[S_{1}\bigtriangleup S_{2}]$ has an induced perfect matching.
\end{theorem}

\setlength {\parindent}{0.0cm}\textbf{Proof.} ($\mathit{i}$) $\Rightarrow $ (%
$\mathit{ii}$), ($\mathit{iii}$) By Proposition \ref{prop2} we get that $%
\Omega (G^{2})\subseteq \Omega (G)$ holds for $G$ square-stable. Now, if $%
S\in \Omega (G^{2})$, and $A$ is a stable set in $G$ disjoint from $S$,
Proposition \ref{prop1} implies that $A$ can be matched into $S$. \mathstrut
If there exists another matching of $A$ into $S$, then at least one vertex $%
a\in A$ has two neighbors in $S$, say $b,c$. Hence, $bc\in E(G^{2})$ and
this contradicts the stability of $S$. Therefore, any $S\in \Omega
(G^{2})\subseteq \Omega (G)$ has property $P1$.%
\setlength
{\parindent}{3.45ex}

($\mathit{ii}$) $\Rightarrow $ ($\mathit{i}$) Suppose, on the contrary, that 
$G$ is not square-stable. It follows that $S_{0}\notin \Omega (G^{2})$,
i.e., there are $v,w\in S_{0}$ with $vw\in E(G^{2})$. Henceforth, there
exists $u\in V-\{v,w\}$, such that $uv,uw\in E(G)$. Consequently, there are
two matchings of $A=\{u\}$ into $S_{0}$, contradicting the fact that $S_{0}$
has property $P1$.

($\mathit{iii}$) $\Rightarrow $ ($\mathit{iv}$) Let $S_{1}\in \Omega (G)$
and $S_{2}\in \Omega (G^{2})$. Then $|S_{2}|\leq |S_{1}|$, and since $%
S_{1}-S_{2}$ is stable in $G$ and disjoint from $S_{2}$, we infer that $%
S_{1}-S_{2}$ can be uniquely matched into $S_{2}$, precisely into $%
S_{2}-S_{1}$, and because $\left| S_{2}-S_{1}\right| \leq \left|
S_{1}-S_{2}\right| $, this matching is perfect. In conclusion, $%
G[S_{1}\bigtriangleup S_{2}]$ has a unique perfect matching.

($\mathit{iv}$) $\Rightarrow $ ($\mathit{v}$) It is clear.

($\mathit{v}$) $\Rightarrow $ ($\mathit{i}$) If $G[S_{1}\bigtriangleup
S_{2}] $ has a perfect matching, for any $S_{1}\in \Omega (G)$ and $S_{2}\in
\Omega (G^{2})$, it follows that $|S_{1}-S_{2}|=|S_{2}-S_{1}|$, and this
implies $|S_{1}|=|S_{2}|$, i.e., $\alpha (G)=\alpha (G^{2})$ is valid.

($\mathit{i}$) $\Rightarrow $ ($\mathit{vi}$) According to ($\mathit{iv}$), $%
G[S_{1}\bigtriangleup S_{2}]$ has a unique perfect matching $M$, for any $%
S_{1}\in \Omega (G)$ and $S_{2}\in \Omega (G^{2})$. By ($\mathit{ii}$), $%
\left| N(v)\cap S_{2}\right| =1$ holds for any $v\in S_{1}-S_{2}$.
Therefore, $M$ must be induced.

($\mathit{vi}$) $\Rightarrow $ ($\mathit{iv}$) It is evident. \rule{2mm}{2mm}

\begin{corollary}
There are no $\alpha ^{-}$-stable square-stable graphs.
\end{corollary}

\setlength {\parindent}{0.0cm}\textbf{Proof.} According to Proposition \ref
{prop11}, $G$ is $\alpha ^{-}$-stable provided $\left| N(v)\cap S\right|
\geq 2$ holds for every $S\in \Omega (G)$ and any $v\in V(G)-S$. If $G$ is
also square-stable, then there exists some $S_{0}\in \Omega (G)$ satisfying
property $P1$, which implies that $\left| N(v)\cap S_{0}\right| =1$ holds
for any $v\in V(G)-S_{0}$. This incompatibility concerning $S_{0}$ proves
that $G$ can not be simultaneously square-stable and $\alpha ^{-}$-stable. 
\rule{2mm}{2mm}\setlength {\parindent}{3.45ex}\newline

\begin{figure}[h]
\setlength{\unitlength}{1cm}%
\begin{picture}(5,1)\thicklines

  \multiput(4,0)(1,0){3}{\circle*{0.29}}
  \put(5,1){\circle*{0.29}}
  \put(4,0){\line(1,0){2}}
  \put(5,0){\line(0,1){1}} 
  \put(6,0){\line(-1,1){1}} 
  \put(4,0){\line(1,1){1}}
  \multiput(7,0)(1,0){3}{\circle*{0.29}}
  \multiput(7,1)(1,0){3}{\circle*{0.29}}
  \put(7,0){\line(1,0){2}}
  \put(7,1){\line(1,0){2}}
  \put(7,0){\line(0,1){1}}
  \put(9,0){\line(0,1){1}}

 \end{picture}
\caption{$\alpha ^{-}$-stable graphs : diamond and $C_{6}$.}
\label{fig3}
\end{figure}
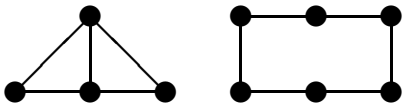

In Figure \ref{fig3} are shown two non-square-stable graphs: $C_{6}$, which
is both $\alpha ^{-}$-stable and $\alpha ^{+}$-stable, and the diamond,
which is only $\alpha ^{-}$-stable.

\begin{corollary}
\label{cor2}Any square-stable graph is $\alpha ^{+}$-stable.
\end{corollary}

\setlength {\parindent}{0.0cm}\textbf{Proof.} Suppose that $G$ is a non-$%
\alpha ^{+}$-stable square-stable graph. Hence, according to Proposition \ref
{prop11}, there are $a,b\in \cap \{S:S\in \Omega (G)\}$, and since $G$ is
square-stable, we infer that $a,b\in \cap \{S:S\in \Omega (G^{2})\}$, as
well. Let $S_{0}\in \Omega (G^{2})$ and $c\in N(a)$ in $G$. Clearly, $a,b\in
S_{0}$, and by Lemma \ref{lem4}, $dist_{G}(a,v)\geq 3$ holds for any $v\in
S_{0}-\{a\}$. Consequently, $dist_{G}(c,v)\geq 2$ holds for any $v\in
S_{0}-\{a\}$. It follows that $S_{1}=S_{0}\cup \{c\}-\{a\}\in \Omega (G)$,
but this contradicts the assumption on $a$, namely that $a\in \cap \{S:S\in
\Omega (G)\}$. \rule{2mm}{2mm}\setlength {\parindent}{3.45ex}\newline

Moreover, we can strengthen Corollary \ref{cor2} to the following:

\begin{corollary}
\label{prop5}Any square-stable graph is well-covered.
\end{corollary}

\setlength {\parindent}{0.0cm}\textbf{Proof.} Assume, on the contrary, that $%
G$ is not well-covered, i.e., there is some maximal stable set $A$ that is
not maximum. According to Theorem \ref{th1}, for any $S\in \Omega (G^{2})$,
there is a unique matching of $B=A-S\cap A$ into $S$, in fact, into $S-A$.
Consequently, $S\cup B-N(B)\cap S$ is a maximum stable set of $G$ that
includes $A$, in contradiction with the fact that $A$ is a maximal stable
set. \rule{2mm}{2mm}\setlength {\parindent}{3.45ex}\newline

It is also possible to see the above result stated implicitly in the proof
of Theorem \ref{th3} from \cite{ranvol}, but our proof is different.

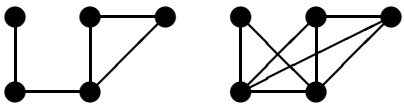
\begin{figure}[h]
\setlength{\unitlength}{1cm}%
\begin{picture}(5,1.2)\thicklines

  \multiput(4,0)(1,0){2}{\circle*{0.29}}
  \multiput(4,1)(1,0){3}{\circle*{0.29}}
  \put(4,0){\line(1,0){1}}
  \put(5,1){\line(1,0){1}}
  \put(5,0){\line(0,1){1}} 
  \put(5,0){\line(1,1){1}} 
  \put(4,0){\line(0,1){1}}

  \multiput(7,0)(1,0){2}{\circle*{0.29}}
  \multiput(7,1)(1,0){3}{\circle*{0.29}}
  \put(7,0){\line(1,0){1}}
  \put(8,1){\line(1,0){1}}
  \put(7,0){\line(1,1){1}} 
  \put(7,0){\line(0,1){1}}
  \put(7,1){\line(1,-1){1}}
  \put(8,0){\line(0,1){1}}
  \put(8,0){\line(1,1){1}}
  \put(7,0){\line(2,1){2}}
 \end{picture}
\caption{A square-stable graph $G$ and its $G^{2}$. $G$ is not very
well-covered.}
\label{fig111111}
\end{figure}

The converse of Corollary \ref{prop5} is not generally true; e.g., $C_{5}$
is well-covered, but is not square-stable. The square-stable graphs do not
coincide with the very well-covered graphs. For instance, $P_{4}$ is both
square-stable and very well-covered, $C_{4}$ is very well-covered and
non-square-stable, but there are square-stable graphs that are not very
well-covered; e.g., the graph in Figure \ref{fig111111}.

\begin{corollary}
Any square-stable graph is $\alpha _{0}^{+}$-stable.
\end{corollary}

\begin{theorem}
\label{th2}For a graph $G$ the following statements are equivalent:

($\mathit{i}$) $G$ is square-stable;

($\mathit{ii}$) there is $S_{0}\in \Omega (G)$ that has the property

$P2:$ \textit{for any stable set} $A$ \textit{of} $G$ \textit{disjoint from} 
$S_{0},$ $A\cup S^{*}\in \Omega (G)$ \textit{holds for some} $S^{*}\subset
S_{0}$;

($\mathit{iii}$) every $S\in \Omega (G^{2})$ has property $P2$.
\end{theorem}

\setlength {\parindent}{0.0cm}\textbf{Proof.} ($\mathit{i}$) $\Rightarrow $ (%
$\mathit{ii}$), ($\mathit{iii}$) By Theorem \ref{th1}, for every $S\in
\Omega (G^{2})$ and any stable set $A$ in $G$, disjoint from $S$, there is a
unique matching of $A$ into $S$. Consequently, $S^{*}=S-N(A)\cap S$ has $%
|S^{*}|=|S|-|A|$ and $S^{*}\cup A\in \Omega (G)$.%
\setlength
{\parindent}{3.45ex}

($\mathit{ii}$) $\Rightarrow $ ($\mathit{i}$) It suffices to show that $%
S_{0}\in \Omega (G^{2})$. If $S_{0}\notin \Omega (G)$, there must exist $%
a,b\in S_{0}$ such that $ab\in E(G^{2})$, and this is possible provided $%
a,b\in N(c)\cap S_{0}$ for some $c\in V-S_{0}$. Hence, $|S_{0}\cup
\{c\}-\{a,b\}|<|S_{0}|$ and this implies that $\{c\}\cup S^{*}\notin \Omega
(G)$ holds for any $S^{*}\subset S$, contradicting the fact that $S_{0}$ has
the property $P2$. Therefore, we get that $S_{0}\in \Omega (G^{2})$, and
this implies that $\alpha (G)=\alpha (G^{2})$.

($\mathit{iii}$) $\Rightarrow $ ($\mathit{i}$) Let $S\in \Omega (G^{2})$, $%
b\in S$ and $a\in V-S$ be such that $ab$ is an edge in $G$. Since $\{a\}$ is
stable and disjoint from $S$, and $S$ has property $P2$, there exists $%
S^{*}\subset S$ so that $S^{*}\cup \{a\}\in \Omega (G)$. Hence, $\left|
S^{*}\right| =\alpha (G)-1$ and consequently, $\left| S\right| =\left|
S^{*}\right| +1=\alpha (G)$, i.e., $S\in \Omega (G)$ holds for any $S\in
\Omega (G^{2})$. By Proposition \ref{prop2}, $G$ is square-stable. \rule%
{2mm}{2mm}\setlength {\parindent}{3.45ex}\newline

Combining Theorem \ref{th3} and our results on square-stable graphs, we
obtain:

\begin{theorem}
\label{th5}For a graph $G$ the following statements are equivalent:

($\mathit{i}$) every vertex of $G$ belongs to exactly one simplex of $G;$

($\mathit{ii}$) $G$ is square-stable;

($\mathit{iii}$) $G$ satisfies $\theta (G)=\theta (G^{2})$;

($\mathit{iv}$) $G$ satisfies $\alpha (G^{2})=\theta (G^{2})=\gamma
(G)=i(G)=\alpha (G)=\theta (G)$;

($\mathit{v}$) $\Omega (G^{2})\subseteq \Omega (G)$;

($\mathit{vi}$) there is some $S\in \Omega (G)$ such that $dist(a,b)\geq 3$
holds for any distinct $a,b\in S$;

($\mathit{vii}$) there exists $S_{0}\in \Omega (G)$ that satisfies the
property $P1$;

$P1$: \textit{any stable set }$A$ \textit{of} $G$ \textit{disjoint from }$%
S_{0}$ \textit{can be uniquely matched into }$S_{0}$;

($\mathit{viii}$) any $S\in \Omega (G^{2})$ has property $P1$;

($\mathit{ix}$) for any $S_{1}\in \Omega (G)$ and $S_{2}\in \Omega (G^{2})$, 
$G[S_{1}\bigtriangleup S_{2}]$ has a unique\emph{\ }perfect matching;

($\mathit{x}$) for any $S_{1}\in \Omega (G)$ and $S_{2}\in \Omega (G^{2})$, $%
G[S_{1}\bigtriangleup S_{2}]$ has a\emph{\ }perfect matching;

($\mathit{xi}$) for any $S_{1}\in \Omega (G)$ and $S_{2}\in \Omega (G^{2})$, 
$G[S_{1}\bigtriangleup S_{2}]$ has an induced perfect matching;

($\mathit{xii}$) there is $S_{0}\in \Omega (G)$ that has the property

$P2:$ \textit{for any stable set }$A$\textit{\ of }$G$\textit{\ disjoint
from }$S_{0}$\textit{, }$A\cup S^{*}\in \Omega (G)$\textit{\ holds for some }%
$S^{*}\subset S_{0}$;

($\mathit{xiii}$) any $S\in \Omega (G^{2})$ has property $P2$.
\end{theorem}

We can now characterize the square-stable graphs that are also simplicial or
chordal, by extending two results from \cite{pristopvest}.

\begin{proposition}
For a graph $G$ the following assertions are equivalent:

($\mathit{i}$) $G$ is square-stable;

($\mathit{ii}$) $G$ is simplicial and well-covered;

($\mathit{iii}$) every vertex belongs to exactly one simplex of $G$.
\end{proposition}

\setlength {\parindent}{0.0cm}\textbf{Proof.} The equivalence ($\mathit{ii}$%
) $\Leftrightarrow $ ($\mathit{iii}$) is proved in \cite{pristopvest}, and
Theorem \ref{th5} ensures that ($\mathit{i}$) $\Leftrightarrow $($\mathit{iii%
}$). \rule{2mm}{2mm}\setlength {\parindent}{3.45ex}

\begin{proposition}
For a chordal graph $G$ the following assertions are equivalent:

($\mathit{i}$) $G$ is square-stable;

($\mathit{ii}$) $G$ is well-covered;

($\mathit{iii}$) every vertex belongs to exactly one simplex of $G$.
\end{proposition}

\setlength {\parindent}{0.0cm}\textbf{Proof.} The equivalence $(\mathit{ii}%
)\Leftrightarrow (\mathit{iii})$ is proved in \cite{pristopvest}, and
Theorem \ref{th5} ensures that ($\mathit{i}$) $\Leftrightarrow $ ($\mathit{%
iii}$). \rule{2mm}{2mm}\setlength {\parindent}{3.45ex}\newline

As another consequence of Theorem \ref{th5}, we obtain that $\Omega (G)$ is
the set of bases of a matroid on $V(G)$ provided $G$ is a complete graph.

\begin{lemma}
\label{lem2}$\Omega (G)$ is the set of bases of a matroid on $V$ if and only
if $\Omega (G^{2})=\Omega (G)$.
\end{lemma}

\setlength {\parindent}{0.0cm}\textbf{Proof.} If $\Omega (G)$ is the set of
bases of a matroid on $V$, then any $S\in \Omega (G)$ must have the property 
$P2$. By Theorem \ref{th2}, $G$ is square-stable and therefore $\Omega
(G^{2})\subseteq \Omega (G)$. Suppose that there exists $S_{0}\in \Omega
(G)-\Omega (G^{2})$; it follows that there are $a,b\in S_{0}$ and $c\in
N(a)\cap N(b)$. Hence, $\{c\}$ is stable in $G$ and disjoint from $S_{0}$,
but $S^{*}\cup \{c\}\notin \Omega (G)$ for any $S^{*}\subset S_{0}$, a
contradiction, since $S_{0}$ has property $P2$. Consequently, the equality $%
\Omega (G^{2})=\Omega (G)$ is true.\setlength {\parindent}{3.45ex}

Conversely, according to Theorem \ref{th2}, any $S\in \Omega (G^{2})=\Omega
(G)$ has the property $P2$. Therefore, $\Omega (G)$ is the set of bases of a
matroid on $V$. \rule{2mm}{2mm}\newline

Combining Theorem \ref{th4} and Lemma \ref{lem2}, we get the following:

\begin{proposition}
\label{prop4}$\Omega (G)$ is the set of bases of a matroid on $V$ if and
only if $G$ is a complete graph.
\end{proposition}

For graphs that are not necessarily connected, we may deduce the following:

\begin{proposition}
\cite{ding} $\Omega (G)$ is the set of bases of a matroid on $V(G)$ if and
only if $G$ is a disjoint union of cliques.
\end{proposition}

\section{Unique pendant perfect matching graphs}

In general, a graph having a unique perfect matching is not necessarily
square-stable. For instance, $K_{3}+e$ has a unique perfect matching, but is
not square-stable. Further, we pay attention to graphs having a perfect
matching consisting of pendant edges, which is obviously unique.

\begin{proposition}
\label{prop6}If $G$ has a perfect matching consisting of pendant edges, then
the following statements are valid:

($\mathit{i}$) $\Omega (G^{2})=\{S_{0}\}$, where $S_{0}=\{v:v$ \textit{is a
pendant vertex in }$G\}$;

($\mathit{ii}$) $G$ is square-stable.
\end{proposition}

\setlength {\parindent}{0.0cm}\textbf{Proof.} Clearly, the set $S_{0}=\{v:v$ 
\textit{is a pendant vertex in }$G\}$ is stable in $G$, and $\left|
S_{0}\right| =\left| V-S_{0}\right| \leq \alpha (G)$. Let $S_{1}\in \Omega
(G)$ and suppose that $\left| S_{1}\right| >\left| S_{0}\right| $. Hence,
both $S_{1}\cap S_{0}$ and $S_{1}\cap (V-S_{0})$ are non-empty, and $\left|
S_{1}\cap S_{0}\right| >\left| V-S_{0}-(S_{1}\cap (V-S_{0}))\right| $. In
addition, we have that $(S_{1}\cap S_{0},S_{1}\cap (V-S_{0}))=\emptyset $,
because $S_{1}$ is stable, and therefore $S_{1}\cap S_{0}$ can not be
matched into $V-S_{0}-(S_{1}\cap (V-S_{0}))$, contradicting the fact that $G$
has a perfect matching. Consequently, $S_{0}\in \Omega (G)$, and because $%
dist_{G}(a,b)\geq 3$ holds for any $a,b\in S_{0}$, we get that $S_{0}\in
\Omega (G^{2})$, i.e., $G$ is square-stable.\setlength {\parindent}{3.45ex}

Assume that there is $S_{2}\in \Omega (G^{2}),S_{0}\neq S_{2}$. Then $%
S_{2}\in \Omega (G)$ and $dist_{G}(a,b)\geq 3$ holds for any $a,b\in S_{2}$.
Let denote $S_{0}=\{v_{i}:1\leq i\leq \alpha (G)\}$ and $N(v_{i})=\{w_{i}\}$%
, for $1\leq i\leq \alpha (G)$. Since $S_{0}\neq S_{2}$, we may assume that,
for instance, $w_{1}\in S_{2}$, and because $w_{1}$ is not pendant, it
follows that $\left| N(w_{1})\right| \geq 2$. Without loss of generality, we
may suppose that $w_{2}\in N(w_{1})$. Hence, $v_{1},v_{2},w_{2}\notin S_{2}$%
, and this implies that $\left| S_{2}\right| <\left| S_{0}\right| $, because
for any $i\geq 3$, $S_{2}$ contains either $v_{i}$ or $w_{i}$, but never
both of them. So, we may conclude that $\Omega (G^{2})=\{S_{0}\}$. \rule%
{2mm}{2mm}\newline

Let us notice that there are square-stable graphs with more than one maximum
stable set, and having no perfect matching; e.g., the graph in Figure \ref
{fig111111}.

\begin{proposition}
\label{prop9}For a K\"{o}nig-Egerv\'{a}ry graph $G$ of order $n\geq 2$ the
following assertions are equivalent:

($\mathit{i}$) $G$ square-stable;

($\mathit{ii}$) $G$ has a perfect matching consisting of pendant edges;

($\mathit{iii}$) $G$ is very well-covered with exactly $\alpha (G)$ pendant
vertices.
\end{proposition}

\setlength {\parindent}{0.0cm}\textbf{Proof.} ($\mathit{i}$) $\Rightarrow $ (%
$\mathit{ii}$) By Proposition \ref{prop5}, $G$ is well-covered, and
according to Proposition \ref{prop8} it is also very well-covered. Hence, we
get that $\alpha (G)=\mu (G)=n/2$, and $G$ has a perfect matching $M$. Let $%
S_{0}=\{a_{i}:1\leq i\leq \alpha (G)\}\in \Omega (G^{2})$ and $b_{i}\in
V(G)-S_{0}$ be such that $a_{i}b_{i}\in M$ for $1\leq i\leq \alpha (G)$. By
Proposition \ref{prop3}, $dist_{G}(v,w)\geq 3$ holds for any $v,w\in S_{0}$.
We claim that every $a_{i}\in S_{0}$ is pendant, i.e., $N(a_{i})=\{b_{i}\}$,
since otherwise, if $b_{j}\in N(a_{i})$ for some $i\neq j$, it follows that $%
dist_{G}(a_{i},a_{j})=2$, in contradiction with $dist_{G}(a_{i},a_{j})\geq 3$%
. Therefore, $M$ consists only of pendant edges.%
\setlength
{\parindent}{3.45ex}

($\mathit{ii}$) $\Rightarrow $ ($\mathit{iii}$) Let $M=\{v_{i}w_{i}:1\leq
i\leq n/2\}$ be the perfect matching of $G$, consisting only of pendant
edges, and suppose that all vertices in $S_{0}=\{v_{i}:1\leq i\leq n/2\}$
are pendant. By Proposition \ref{prop6}, we get that $S_{0}\in \Omega (G)$,
i.e., $\alpha (G)=\mu (G)=n/2$. Assume that $G$ is not well-covered, that is
there exists some maximal stable set $A$ in $G$ such that $A\notin \Omega
(G) $. Since $S_{0}$ contains all pendant vertices of $G$, it follows that $%
A\cup \{v_{i}:v_{i}\in S_{0},N(v_{i})\cap A=\emptyset \}$ is stable and
larger than $A$, in contradiction with the maximality of $A$. In conclusion, 
$G$ is very well-covered.

($\mathit{iii}$) $\Rightarrow $ ($\mathit{i}$) Since $G$ is very
well-covered with exactly $\alpha (G)$ pendant vertices, we infer that $%
S_{0}=\{v:v$ \textit{is a pendant vertex}$\}\in \Omega (G)$ and also that
the matching $M=\{vw:vw\in E(G),v\in S_{0}\}$ is perfect and consists of
only pendant edges. According to Proposition \ref{prop6}, it follows that $G$
is square-stable. \rule{2mm}{2mm}

\begin{remark}
Well covered K\"{o}nig-Egerv\'{a}ry graphs do not have to be square-stable,
for instance, the graph $C_{4}$.
\end{remark}

\begin{remark}
A K\"{o}nig-Egerv\'{a}ry graph with a unique perfect matching is not always
square-stable, e.g., the graphs $P_{6}$ (by the way, it is also a tree) and $%
K_{3}+e$ (i.e., the graph in Figure \ref{fig2}).
\end{remark}

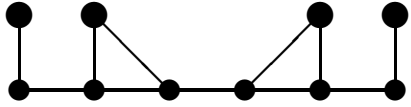
\begin{figure}[h]
\setlength{\unitlength}{1.0cm} 
\begin{picture}(5,1)\thicklines
  \multiput(4,0)(1,0){6}{\circle*{0.29}}
  \multiput(4,1)(1,0){2}{\circle*{0.35}}
  \multiput(8,1)(1,0){2}{\circle*{0.35}}
  \put(4,0){\line(1,0){5}}
  \put(5,1){\line(1,-1){1}}
  \put(7,0){\line(1,1){1}}
  \multiput(4,0)(1,0){2}{\line(0,1){1}}
  \multiput(8,0)(1,0){2}{\line(0,1){1}}
\end{picture}
\caption{$G$ is square-stable and has a unique perfect matching containing
not only pendant edges.}
\label{fig12222}
\end{figure}

\begin{remark}
A non-K\"{o}nig-Egerv\'{a}ry graph with a unique perfect matching $M$ may be
square-stable, even if $M$ does not consist of only pendant edges (for
instance, see the graph in Figure \ref{fig12222}).
\end{remark}

Proposition \ref{prop9} is true for bipartite graphs as well, since any
bipartite graph is also a K\"{o}nig-Egerv\'{a}ry graph. It is worth
recalling here that for a bipartite graph (see \cite{levm1}, and for trees
see \cite{gun}) to have a perfect matching is equivalent to be $\alpha ^{+}$%
-stable. In general, we have shown in \cite{levm2} that any $\alpha ^{+}$%
-stable K\"{o}nig-Egerv\'{a}ry graph has a perfect matching, while the
converse is not true (see, for instance, the diamond, Figure \ref{fig3}).

\begin{proposition}
\cite{levm3}\label{prop10} Any well-covered tree $T$ non-isomorphic to $%
K_{1},K_{2}$, contains at least one edge $e$ connecting two non-pendant
vertices, such that $T-e=T^{\prime }\cup K_{2}$ and $T^{\prime }$ is a
well-covered tree.
\end{proposition}

For trees, Propositions \ref{prop9} and \ref{prop10} lead to the following
extension of the characterization that Ravindra gave to well-covered trees
in \cite{rav1}:

\begin{corollary}
\label{cor1}If $T$ is a tree of order $n\geq 2$, then the following
statements are equivalent:

($\mathit{i}$) $T$ is well-covered;

($\mathit{ii}$) $T$ is very well-covered;

($\mathit{iii}$) $T$ has a perfect matching consisting of pendant edges;

($\mathit{iv}$) $T$ is square-stable.
\end{corollary}

\setlength {\parindent}{0.0cm}\textbf{Proof.} Let us notice that for general
graphs: ($\mathit{iv}$) $\Rightarrow $ ($\mathit{i}$) is true according to
Corollary \ref{prop5}, and the implication ($\mathit{iii}$) $\Rightarrow $ ($%
\mathit{ii}$) is clear. Further, for K\"{o}nig-Egerv\'{a}ry graphs, the
assertions ($\mathit{iii}$), ($\mathit{iv}$) are equivalent according to
Proposition \ref{prop9}, and ($\mathit{i}$), ($\mathit{ii}$) are equivalent
by Proposition \ref{prop12}. Thus, to complete the proof of the corollary,
it is sufficient to show that for trees ($\mathit{i}$) implies ($\mathit{iii}
$). Since ($\mathit{i}$) and ($\mathit{ii}$) are equivalent, the order $n$
of $T$ must be even. We use induction on $n$. The assertion is true for $n=2$%
. If $T$ has $n>2$ vertices, then according to Proposition \ref{prop10}, $T$
contains at least one edge $e$ connecting two non-pendant vertices, such
that $T-e=T^{\prime }\cup K_{2}$ and $T^{\prime }$ is a well-covered tree.
By the induction hypothesis, $T^{\prime }$ has a perfect matching $M$
consisting of pendant edges. Hence, $M\cup \{e\}$ is a perfect matching of $%
T $ consisting of pendant edges. \rule{2mm}{2mm}%
\setlength
{\parindent}{3.45ex}\newline

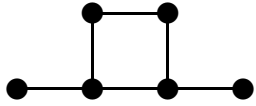
\begin{figure}[h]
\setlength{\unitlength}{1cm}%
\begin{picture}(5,1)\thicklines

  \multiput(5,0)(1,0){4}{\circle*{0.29}}
  \multiput(6,1)(1,0){2}{\circle*{0.29}}
  \put(5,0){\line(1,0){3}}
   
  \put(6,0){\line(0,1){1}} 
  \put(7,0){\line(0,1){1}} 
  \put(6,1){\line(1,0){1}} 
  
 \end{picture}
\caption{A bipartite and very well-covered but not square-stable graph.}
\label{fig4}
\end{figure}

Let us notice that the equivalences appearing in Corollary \ref{cor1} fail
for bipartite graphs. For instance, the graph in Figure \ref{fig4}\emph{\ }%
is very well-covered, but is not square-stable. 

Combining Proposition \ref{prop8} and Proposition \ref{prop9}, we obtain:

\begin{corollary}
$G$ is square-stable and very well-covered if and only if $G$ is a
K\"{o}nig-Egerv\'{a}ry graph with exactly $\alpha (G)$ pendant vertices.
\end{corollary}

\begin{corollary}
If $G$ is a square-stable K\"{o}nig-Egerv\'{a}ry graph, then $G^{2}$ is also
a K\"{o}nig-Egerv\'{a}ry graph.
\end{corollary}

\begin{remark}
Figure \ref{fig4} brings an example of a K\"{o}nig-Egerv\'{a}ry graph whose
square is not a K\"{o}nig-Egerv\'{a}ry graph.
\end{remark}

Another consequence of Proposition \ref{prop9} is the following extension of
the characterization that Finbow, Hartnell and Nowakowski give in \cite
{finhart1} for graphs having the girth $\geq 6$.

\begin{proposition}
Let $G$ be a graph of girth $\geq 6$, which is isomorphic to neither $C_{7}$
nor $K_{1}$. Then the following assertions are equivalent:

($\mathit{i}$) $G$ is well-covered;

($\mathit{ii}$) $G$ has a perfect matching consisting of pendant edges;

($\mathit{iii}$) $G$ is very well-covered;

($\mathit{iv}$) $G$ is a K\"{o}nig-Egerv\'{a}ry graph with exactly $\alpha
(G)$ pendant vertices;

($\mathit{v}$) $G$ is a K\"{o}nig-Egerv\'{a}ry square-stable graph.
\end{proposition}

\setlength {\parindent}{0.0cm}\textbf{Proof.} The equivalences ($\mathit{i}$%
) $\Leftrightarrow $ ($\mathit{ii}$) $\Leftrightarrow $ ($\mathit{iii}$) are
done in \cite{finhart1}. In \cite{levm} it has been proved that ($\mathit{iii%
}$) $\Leftrightarrow $ ($\mathit{iv}$). Finally, ($\mathit{ii}$) $%
\Leftrightarrow $ ($\mathit{v}$) is true by Propositions \ref{prop6} and \ref
{prop9}. \rule{2mm}{2mm}\setlength {\parindent}{3.45ex}

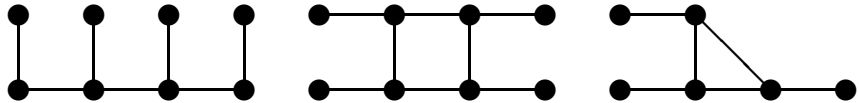
\begin{figure}[h]
\setlength{\unitlength}{1.0cm} 
\begin{picture}(5,1.2)\thicklines
  \multiput(1,0)(1,0){12}{\circle*{0.29}}
  \multiput(1,1)(1,0){10}{\circle*{0.29}}
  \put(1,0){\line(1,0){3}}
  \put(5,0){\line(1,0){3}}
  \put(5,1){\line(1,0){3}}
  \put(9,0){\line(1,0){3}}
  \put(9,1){\line(1,0){1}}
  \multiput(1,0)(1,0){4}{\line(0,1){1}}
  \multiput(6,0)(1,0){2}{\line(0,1){1}}
  \put(10,0){\line(0,1){1}}
  \put(10,1){\line(1,-1){1}} 
\end{picture}
\caption{Square-stable Koenig-Egervary graphs.}
\label{fig1212}
\end{figure}

\begin{remark}
$C_{7}$ is not a K\"{o}nig-Egerv\'{a}ry graph.
\end{remark}

\section{Conclusions}

In this paper we continue the investigations, started by Randerath and
Volkmann \cite{ranvol}, on the class of square-stable graphs. We think that
the characterization of Koenig-Egervary square-stable graphs obtained here
may be extended to some new classes of square-stable graphs. It is also
important to mention that square-stable trees have a very specific recursive
structure (see \cite{levm3}).

It also seems interesting to study graphs satisfying some equalities between
the invariants appearing in the following series of inequalities: $\alpha
(G^{2})\leq \theta (G^{2})\leq \gamma (G)\leq i(G)\leq \alpha (G)\leq \theta
(G)$, for instance $\alpha (G^{2})=i(G)$.

\end{document}